\documentclass[a4paper, 12pt]{article}
\usepackage{enumerate, theorem}
\usepackage{amsmath, amsfonts, amssymb}
\usepackage[height=22.5cm, width=15cm]{geometry}
\usepackage[all]{xy}
\usepackage{mathrsfs, yfonts}

\DeclareMathOperator{\id}{id}

 \def\C{{\mathbb C}}

\newcommand{\parag}[1]{\paragraph{\sc{#1.}} }
\newcommand{\A}{\tilde{\mathcal{A}}}
\newcommand{\Ab}{\mathcal{A}}

\newtheorem{thm}{Th\'eor\`{e}me}[subsection]

\newtheorem{prop}[thm]{Proposition}
\newtheorem{lemma}[thm]{Lemme}

\begin{document}

\date{24/02/14}

\author{Daniel Barlet\footnote{Barlet Daniel, Institut Elie Cartan UMR 7502  \newline
Universit\'e de Lorraine, CNRS, INRIA  et  Institut Universitaire de France, \newline
BP 239 - F - 54506 Vandoeuvre-l\`es-Nancy Cedex.France. \newline
e-mail : daniel.barlet@univ-lorraine.fr}.}

\title{ Complement  for Algebraic differential equations...}

\maketitle

\tableofcontents

\parag{Abstract} We complete the study of some periods of polynomials in (n+1)-variables with (n+2)-monomials in computing the behavior of these periods in the natural parameter for such a polynomial.

\parag{Classification AMS 2010} 32-S-30, 34-M-56, 58-K-60.

\newpage

\section{Introduction}

This note is a complement to the study in [B.13]  of period integrals of non quasi-homogeneous polynomials in \ $n+1$ \ variables with \ $n+2$ \ monomials. We focuse here on the dependance of these period integrals on the ``natural'' parameter \ $\lambda \in \C^{*}$ \ which is the only ``free''  coefficient of such a polynomial modulo the dilatations of the variables\footnote{In our hypothesis we may assume that all monomials have coefficient \ $1$ \ excepted the last one up to a linear diagonal change of variable.}\\
For that purpose we recall first the fact that for a polynomial function \ $f$ \ depending polynomially of a parameter \ $\lambda$ \ we may define a natural \ ``$b-$connection'' on the highest \ $(f,\lambda)-$relative de Rham cohomology group of \ $f$ \ which induces the derivation \ $\frac{\partial }{\partial \lambda} $ \ on period integrals. The construction for any holomorphic function depending of a holomorphic parameter is precised in the appendix.\\
Then we show how to  compute explicitely this connection in our specific situation and we obtain a simple partial differential equation for the period integrals associated to any monomial in \ $\C[x_{0},\dots, x_{n}]$ \ when we consider a polynomial of the type
$$ f_{\lambda}(x) = \sum_{j = 1}^{n+1} \  x^{\alpha_{j}} \ + \lambda.x^{\alpha_{n+2}} \qquad {\rm where} \qquad \alpha_{j} \in \mathbb{N}^{n+1}, \ j \in [1,n+2]$$
with he following assumptions
\begin{enumerate}[i)]
\item The \ $(n+2,n+2)-$matrix obtained from \ $M : = (\alpha_{1}, \dots, \alpha_{n+2})$ \ by adding a first line of \ $1$ \ has rank \ $n+2$.
\item The elements \ $\alpha_{1}, \dots, \alpha_{n+1}$ \ form a basis of \ $\mathbb{Q}^{n+1}$.
\end {enumerate}
Note that the first condition is equivalent to the fact that \ $f$ \ is not quasi-homogeneous, and that the condition ii) is always satisfied assuming i), up to change the order of the monomials (and then to change the parameter \ $\lambda$ \ to \ $c.\lambda^{m}$ \ for some \ $c \in \mathbb{Q}^{*}$ \ some \ $m \in \mathbb{Z}^{*})$.

\section{The \ $\lambda-$connection.}

\subsection{The general situation.}

We consider here a polynomial \ $f\in R : = \C[x_{0},\dots,x_{n}][\lambda]$ \ depending polynomially on a parameter \ $\lambda$. We consider on \ $R\otimes \Lambda^{*}(\C^{n+1}) : = \Omega_{/}^{*}$ \  the \ $\lambda-$relative de Rham complex, where \ $(\C^{n+1})^{*} : = \oplus_{i=0}^{n} \  \C.dx_{i}$,  and we denote  \ $d_{/}$ \ its differential.\\
We shall denote by \ $\Ab$ \ the unitary (non commutative) algebra generated by \ $a$ \ and \ $b$ \ with the commutation relation \ $a.b - b.a = b^{2}$ \ and by \ $\Ab[\lambda] : = \Ab\otimes_{\C} \C[\lambda]$ \ with its natural structure of algebra  for which \ $\lambda$ \ commutes with \ $a$ \ and \ $b$.\\
Then the quotient \ $E_{f} : = \Omega_{/}^{n+1}\big/d_{/}f\wedge d_{/}\Omega_{/}^{n-1}$ \ has a natural left \ $\Ab[\lambda]-$module structure defined by
\begin{itemize}
\item The action of \ $a$ \ is given by the multiplication by \ $f$.
\item The action of \ $b$ \ is given by \ $d_{/}f \wedge d_{/}^{-1}$.\\
\end{itemize}

Remark that for fixed \ $\lambda$, assuming that \ $f_{\lambda}$ \ has an isolated singularity at \ $0$, the \ $b-$completion \ $E_{f}\otimes_{\C[b]}\C[[b]]$ \ is the usual (formal) Brieskorn module associated to \ $f_{\lambda}$ \ at \ $0$. For a given monomial \ $\mu \in \C[x_{0}, \dots, x_{n}]$ \ the decomposition theorem of  [B.13] (theorem 3.1.2) applies to the the quotient \ $\A\big/\A.P(\mu)$ \ where \ $\A$ \ is the \ $b-$completion of \ $\Ab$ \ and where \ $P(\mu)$ \ is the element in \ $\Ab$ \ constructed in the theorem 1.2.1 of [B.13]. Then  \ $P_{d}(\mu)$ \ is a (left-)multiple of the Bernstein element of \ $\A.[\mu.dx]$ \ in \ $E_{f}\otimes_{\C[b]}\C[[b]]$ \ and it determines a finite set of possible eigenvalues for the monodromies around \ $s = 0$ \ for the period integrals ($\lambda$ \ fixed)
$$ \varphi_{\lambda}(s) = \int_{\gamma_{\lambda,s}} \ \frac{\mu.dx}{d_{/}f} $$
for any horizontal family\ $\gamma_{\lambda,s}$ \ of compact \ $n-$dimensional cycles in the fibers of \ $f_{\lambda}$.\\

It is important to remark that if \ $f_{\lambda}$ \ has a non isolated singularity at the origin, despite the fact that there is no finiteness for the \ $\C[[b]]-$module \ $E_{f}\otimes_{\C[b]}\C[[b]]$, the conclusion above is still valid because the quotient \ $\A\big/\A.P(\mu)$, and so its image in \ $E_{f}\otimes_{\C[b]}\C[[b]]$, is a finite type \ $\C[[b]]-$module. Then the product decomposition \ $P_{d} = (a -r_{1}.b)\dots (a -r_{d}.b)$, where \ $r_{1}, \dots, r_{d}$ \ are (explicitely computable)  rational numbers, gives that the set \ $\{e^{2i\pi.r_{1}}, \dots, e^{2i\pi.r_{d}}\}$ \ contains the spectrum of these monodromies (counting multiplicities).

 \parag{Question} Is it true that \ $P_{d} $ \ is  equal to the Bernstein element\footnote{For a (a,b)-module \ $E$ \ with one generator as a \ $\A-$module, the relation between its Bernstein element \ $P_{d} \in \Ab$ \ and its Bernstein polynomial \ $B$ \  is given by the formula (see [B.09])
 $$ (-b)^{d}.B(-b^{-1}.a) = P_{d} \qquad {\rm where} \quad d \quad {\rm is \ the \ rank\ of} \quad E .$$ }   of the Brieskorn module \ $\A.[\mu.dx]$ \ in \ $E_{f}\otimes_{\C[b]}\C[[b]]$ \ when \ $f_{\lambda}$ \ has an isolated singularity at the origin  ? $\hfill \square$\\

\begin{prop}\label{connexion}
There exists a \ $\C-$linear operator \ $\nabla : E_{f} \to E_{f}$ \ with the following properties :
\begin{enumerate}
\item For \ $\omega = d_{/}\xi \in \Omega_{/}^{n+1}$ \ we have \ $\nabla([\omega]) = [d_{/}f\wedge\frac{\partial \xi}{\partial \lambda} -\frac{\partial f}{\partial \lambda}.\omega] .$
\item The map \ $b^{-1}.\nabla$ \ well defined on \ $ b.\tilde{E}_{f}$ \ where \ $\tilde{E}_{f} : = E_{f}\big/(b-torsion)$, with value in \ $\tilde{E}_{f }$, commutes with \ $a$ \ and \ $b$ \ and is a \ $\lambda-$connection.
\item If \ $(\gamma_{s,\lambda})_{(s,\lambda) \in S\times \Omega}$ \ is a horizontal family of compact \ $n-$cycles in the fibers of \ $(f,\lambda)$ \ over an open set in \ $\C\times \C^{*} \setminus C(f,\lambda)$ \ where \ $C(f,\lambda)$ \ is the set of critical values of the map \ $(f,\id) : \C^{n+1}\times \C^{*} \to \C \times \C^{*}$, we have for any \ $\omega \in \Omega_{/}^{n+1}$ \ the equality
$$  \frac{\partial }{\partial s} \frac{\partial }{\partial \lambda}\big[ \int_{\gamma_{s,\lambda}} \ \frac{\omega}{d_{/}f } \big] =  \int_{\gamma_{s,\lambda}}  \ \frac{ \nabla(\omega)}{d_{/}f }.$$
\end{enumerate}
\end{prop}

\parag{Proof} Remark first that \ $\nabla$ \ is well defined because for \ $\xi = d_{/}\eta$ \ we have 
$$\nabla(d_{/}\xi) = d_{/}f\wedge \frac{\partial (d_{/}\eta)}{\partial \lambda} = d_{/}f\wedge d_{/}(\frac{\partial \eta}{\partial \lambda})$$ 
 so it induces \ $0$ \ in \ $E_{f}$.\\
Let \ $\omega = d_{/}\xi \in  \Omega_{/}^{n+1}$ \ and let \ $d_{/}\eta = d_{/}f \wedge \xi$. Then we have
\begin{align*}
&  \nabla(b.[\omega]) = \nabla(d_{/}\eta) = d_{/}f\wedge\frac{\partial \eta}{\partial \lambda} - \frac{\partial f}{\partial \lambda}.d_{/}\eta \\
& \qquad  = d_{/}f\wedge (\frac{\partial \eta}{\partial \lambda} - \frac{\partial f}{\partial \lambda}.\xi) = b\big[d_{/}(\frac{\partial \eta}{\partial \lambda} - \frac{\partial f}{\partial \lambda}.\xi)\big] \\
& \qquad  =  b\big[\frac{\partial\big(d_{/}f\big)}{\partial \lambda}\wedge \xi + d_{/}f\wedge\frac{\partial \xi}{\partial \lambda} - \frac{\partial f}{\partial \lambda}.d_{/}\xi - d_{/}(\frac{\partial f}{\partial \lambda})\wedge \xi )\big] \\
& \qquad  = b\big[d_{/}(\frac{\partial f}{\partial \lambda}).\xi + \nabla(d_{/}\xi)- d_{/}(\frac{\partial f}{\partial \lambda})\wedge \xi )\big]  = b\big[ \nabla(d_{/}\xi)\big]
\end{align*}
 as \ $d_{/}$ \ and \ $\frac{\partial} {\partial \lambda}$ \ commute. So we have \ $b.\nabla = \nabla.b$.\\
We have also
\begin{align*}
& \nabla(a.[\omega]) = \nabla(f.d_{/}\xi) = \nabla(d_{/}(f.\xi)) - \nabla(d_{/}f\wedge \xi) \\
& \qquad  = d_{/}f\wedge \frac{\partial (f\xi)}{\partial \lambda} - \frac{\partial f}{\partial \lambda}.f.d_{/}\xi - \frac{\partial f}{\partial \lambda}.d_{/}f\wedge \xi - \nabla(b.[\omega]) \\
& \qquad  = a.\nabla([\omega]) - b.\nabla([\omega]).
\end{align*}
This implies the equality \ $a.b^{-1}.\nabla = b^{-1}.\nabla.a$ \ as \ $\C-$linear maps from \ $b.\tilde{E}_{f}$ \ to  \ $\tilde{E}_{f}$. \\

Note that the equalities \ $\nabla.b = b.\nabla $ \ and \ $ \nabla.a = (a - b).\nabla$ \ as \ $\C-$endomorphisms of \ $E_{f}$ \ are more precise than the relations above.\\

Finally let \ $\varphi \in \C[\lambda]$ \ then we have
\begin{align*}
&  \nabla(\varphi.d_{/}\xi) = \nabla(d_{/}(\varphi.\xi)) =  d_{/}f\wedge \frac{\partial \varphi.\xi}{\partial \lambda} - \frac{\partial f}{\partial \lambda}.\varphi.d_{/}\xi \\
& \qquad  =  \frac{\partial \varphi}{\partial \lambda}.(d_{/}f\wedge \xi) + \varphi.\nabla(d_{/}\xi) =  \frac{\partial \varphi}{\partial \lambda}.b[d_{/}\xi]  + \varphi.\nabla(d_{/}\xi)\\
\end{align*}
and this shows that \ $b^{-1}.\nabla$ \ is a \ $\lambda-$connection.\\

Note again that we proved  the equality in \ $E_{f}$ :  \ $\nabla(\varphi.\omega) = \frac{\partial \varphi}{\partial \lambda}.b.\omega + \varphi.\nabla(\omega)$ \ valid for \ $\varphi \in \C[\lambda]$ \ and \ $\omega \in E_{f}$ \  which is more precise.\\

To prove the point \ $3.$ \ of the statement consider \ $\xi \in \Omega_{/}^{n}$ \ and let \ $d$ \ be the total de Rham differential (in \ $x$ \ and \ $\lambda)$. We have
$$ d\xi = d\lambda\wedge \frac{\partial \xi}{\partial \lambda} + d_{/}\xi \quad {\rm and} \quad  df = d\lambda.\frac{\partial f}{\partial \lambda} + d_{/}f .$$
Assume  we can  write \ $d\xi = d\lambda\wedge v + d_{/}f\wedge u $ \ with \ $u, v \in \Omega_{/}^{n}$. Then we obtain
\begin{equation*}
 d\xi = d\lambda\wedge(v - \frac{\partial f}{\partial \lambda}.u) + d_{/}f\wedge u \quad  {\rm with} \quad u = \frac{d_{/}\xi}{d_{/}f} \quad {\rm and} \quad  v = \frac{\partial \xi}{\partial \lambda}.
 \end{equation*}
 If \ $(\gamma_{s,\lambda})$ \ is a horizontal family of compact \ $n-$cycles in the fibers of the map \ $(f,\id) : \C^{n+1}\times \C^{*} \to \C\times \C^{*}$, we shall have
 $$ d\big(\int_{\gamma_{s,\lambda}} \ \xi \big)= \big[\int_{\gamma_{s,\lambda}} (v- \frac{\partial f}{\partial \lambda}.u)\big].d\lambda + \big[\int_{\gamma_{s,\lambda}} \ u\big].ds.$$
 So, has the chain \ $\cup_{s,\lambda} \ \gamma_{s,\lambda}$ \ is proper and without \  $\lambda-$relative boundary we obtain
 \begin{equation*}
   \frac{\partial }{\partial s} \int_{\gamma_{s,\lambda}} \ \xi = \int_{\gamma_{s,\lambda}} \ \frac{d_{/}\xi}{ d_{/}f} \qquad {\rm and} \qquad
  \frac{\partial }{\partial \lambda} \int_{\gamma_{s,\lambda}} \ \xi =  \int_{\gamma_{s,\lambda}} \ \big(\frac{\partial \xi}{\partial \lambda} - \frac{\partial f}{\partial \lambda}. \frac{d_{/}\xi}{d_{/}f}\big) 
 \end{equation*}
 Now consider \ $\omega \in \Omega_{/}^{n+1}$ \ and write \ $\omega = d_{/}\xi$. Then \ $b[\omega] = [d_{/}f \wedge\xi] $ \ and we have
 $$ \int_{\gamma_{s,\lambda}} \ \frac{b[\omega]}{d_{/}f }= \int_{\gamma_{s,\lambda}} \xi \qquad {\rm and} \qquad  \frac{\partial }{\partial s} \int_{\gamma_{s,\lambda}} \  \frac{b[\omega]}{d_{/}f } =   \int_{\gamma_{s,\lambda}} \ \frac{\omega}{d_{/}f }.$$
 So we conclude that we have
 \begin{equation*}
  \frac{\partial }{\partial \lambda} \int_{\gamma_{s,\lambda}} \  \frac{\omega}{d_{/}f} = \frac{\partial }{\partial \lambda} \int_{\gamma_{s,\lambda}}  \xi 
   = \int_{\gamma_{s,\lambda}} \ \big( \frac{\partial \xi}{\partial \lambda} - \frac{\partial f}{\partial \lambda}. \frac{d_{/}\xi}{d_{/}f}\big)  =  \frac{\partial }{\partial \lambda} \int_{\gamma_{s,\lambda}} \ \frac{\nabla[\omega]}{d_{/}f}   \qquad   \qquad \blacksquare
  \end{equation*}

\bigskip

\subsection{The case of a polynomial with  \ $n+2$ \ monomials in \ $n+1$ \ variables.}

So we consider now the case were \ $f : = \sum_{j=1}^{n+2} \ m_{j}$ \ where \ $m_{j} : = x^{\alpha_{j}} \quad  j \in[1,n+1]$ \ and \ $m_{n+2} : = \lambda.x^{\alpha_{n+2}}$ \ with the following hypotheses (see [B. 13]) : the rank of  the square matrix\ $M' : = (\alpha_{1}, \dots, \alpha_{n+1})$ \ is \ $n+1$ \ and the rank of the square matrix \ $\tilde{M}$ \ obtained by adding a first line of \ $1$ \ to the matrix \ $M : = (\alpha_{1}, \dots, \alpha_{n+2})$ \ is \ $n+2$.\\
Recall that if we write (with a minimal positive integer \ $r$) \ $r.\alpha_{n+2} = \sum_{j=1}^{n+1} \ p_{j}.\alpha_{j}$ \ where \ $p_{1}, \dots, p_{n+1}$ \ are in \ $\mathbb{Z}$, and if we define 
 $$d = \inf\{ r - \sum_{j, p_{j} \leq 0} p_{j}, \sum_{j,p_{j} > 0} p_{j}\} \quad {\rm and} \qquad d+h = \sup\{ r - \sum_{j, p_{j} \leq 0} p_{j}, \sum_{j,p_{j} > 0} p_{j}\}$$
  there exists an element \ $P$ \ in \ $A[\lambda,\lambda^{-1}]$  \ of the form
$$ P : = P_{d+h} + c.\lambda^{\pm r}.P_{d}$$
which annihilated the class \ $[dx]$ \ in \ $E_{f}$, where \ $ P_{d+h}$ \ and \ $P_{d}$ \ are homogeneous elements in \ $\Ab$, respectively of degree \ $d+h$ \ and \ $d$ \ which are monic in \ $a$ \ with rational coefficients\footnote{More is proved in [B.13] : \ $P_{d+h}$ \ and \ $P_{d}$ \ factorize in product of \ $(a - r_{j}.b)$ \ with \ $r_{j}\in \mathbb{Q}$.},  and where \ $c$ \ is in \ $\mathbb{Q}^{*}$. The sign in the exponent of \ $\lambda$ \ will be precised in the proof of the proposition \ref{conn. en bas}.\\
Recall also that in this situation the \ $\Ab[\lambda]$ \ module generated by the class \ $[dx]$ \ in \ $E_{f}$ \ is exactely the image in \ $E_{f}$ \ of \ $\C[m_{1},\dots,m_{n+2}][\lambda].dx \subset \Omega_{/}^{n+1}$ \ with \ $m_{j} = x^{\alpha_{j}}$ \ with \ $ j \in [1,n+1]$ \ and \ $m_{n+2} = \lambda.x^{\alpha_{n+2}}$.

\bigskip

Our next result uses the following easy lemma:

\begin{lemma}\label{tres facile}
Let \ $Q \in \Ab$ \ a homogeneous element in \ $(a,b)$ \ of degree \ $k$. Then for any \ $\lambda \in \C$ \ we have :
$$ b.Q.b^{-1}.(a - \lambda.b) = (a -(\lambda+k).b).Q. $$
\end{lemma}

\parag{proof} Remark first that the map \ $\Ab \to \Ab$ \ sending \ $x \in \Ab$ \ to \ $b.x.b^{-1}$ \ is well defined  and bijective thanks to the following facts : \ $b$ \ is injective and \ $b.\Ab = \Ab.b$.\\
We shall prove the lemma by induction \ $k$. As the case \ $k = 0$ \ is obvious, assume that the lemma is proved for \ $k < k_{0}$ \  where \ $k_{0} \geq 1$ \ and consider an homogeneous element \ $Q$ \ of degree \ $k_{0}$. We may assume\footnote{Recall that any homogeneous element in \ $\Ab$ \ which is monic in \ $a$ \ factorizes as a product of linear factors \ $(a - r_{i}.b)$, where the \ $r_{i}$ \ are complex numbers; see [B.09].}  that \ $Q = b.R$ \ or that we may find \ $\mu \in \C$ \ such that \ $Q = (a - \mu.b).R$, where \ $R$ \ is homogeneous of degree \ $k_{0}-1$. In the first case we have, using the induction hypothesis :
$$ b.b.R.b^{-1}.(a - \lambda.b) = b.(a - (\lambda +k_{0}-1).b).R = (a - (\lambda+ k_{0}).b).b.R = (a - (\lambda+ k_{0}).b).Q.$$
In the second case we have, using the induction hypothesis :
\begin{align*}
&  b.(a - \mu.b).R.b^{-1}.(a - \lambda.b)    = (a - (\mu+1).b).b.R.b^{-1}.(a - \lambda.b) \\
& \qquad \qquad \qquad \qquad  \qquad  \quad  = (a -(\mu +1).b).(a - (\lambda+k_{0}-1).R \\
& \qquad \qquad \qquad \qquad \qquad  \quad  = (a -(\lambda+k_{0}).b).(a - \mu.b).R \\
& \qquad \qquad \qquad \qquad \qquad  \quad   =   (a -(\lambda+k_{0}).b).Q.    \hfill \qquad   \qquad  \qquad  \qquad  \qquad  \qquad \qquad  \quad \blacksquare
\end{align*}

\begin{prop}\label{conn. en bas}
Let \ $\mu$ \ be a monomial of degree \ $k$ \ in \ $\C[x_{0}, \dots,x_{n+1}]$. Then we have in \ $E_{f}$ \ the relation
$$ \nabla([\mu]) = \frac{-1}{\lambda}.(\sigma.a + (\tau -k.\sigma).b)[\mu] $$
where \ $\sigma, \tau$ \ are defined by the relation \ $m_{n+2}.[\mu] = (\sigma.a + \tau.b)[\mu]$. Moreover the value\footnote{The sign is precised in the proof and only depends on \ $\alpha_{1}, \dots, \alpha_{n+2}$.} of \ $\sigma$ is \ $\pm r/h$ \ so it does not depend on the choice of the monomial \ $\mu$.\\
As a consequence, if we have on an open set \ $S \times \Omega$ \ in \ $\C^{*}\times \C^{*}$, a horizontal family \ $(\gamma_{s,\lambda})_{(s,\lambda)\in S\times \Omega}$ \ of compact  \ $n-$cycles in the fibers of the map \ $\C^{n+1}\times \C^{*} \to \C\times \C^{*}$ \ defined by \ $(x,\lambda) \mapsto (f_{\lambda}(x), \lambda)$, the holomorphic function 
 $$(s,\lambda) \mapsto \varphi(s, \lambda) : = \int_{\gamma_{s,\lambda}} \ \frac{\mu.dx}{d_{/}f} $$
 satisfies the partial differential equation
 $$- \lambda.\frac{\partial }{\partial \lambda}\frac{\partial }{\partial s} \varphi = \sigma.\frac{\partial (s.\varphi) }{\partial s} + (\tau - k).\varphi $$
 on \ $S \times \Omega$. $\hfill \blacksquare$
\end{prop}

\parag{Proof} As we have \ $\lambda.\nabla([1]) = - m_{n+2} $ \ in \ $E_{f}$ \ and as we know that there exist \ $\sigma,\tau$ \ in \ $\mathbb{Q}$ \ such that \ $ (\sigma.a + \tau.b)[1] = m_{n+2}$ \ for the case \ $\mu = 1$ \ the only thing to prove is the computation of \ $\sigma$.\\
Using the Cramer system with matrix \ $(n+2,n+2)$ \ obtained by adding a first line of \ $1$ \ to the matrix \ $M : = (\alpha_{1}, \dots, \alpha_{n+2})$,  computing \ $a[1]$ \ and the \ $b_{i}[1]$ \ we find that \ $\sigma$ \  is the coefficient \ $(n+2,1)$ \ in the matrix \ $\tilde{M}^{-1}$. Let \ $M'$ \ be the principal \ $(n+1,n+1)$ \ minor of \ $\tilde{M}$. This implies that
$$ \sigma = (-1)^{n+1} \frac{det(M')}{det(\tilde{M})} .$$
But using the relation \ $\alpha_{n+2} = \sum_{j=1}^{n+1} \ \frac{p_{j}}{r}.\alpha_{j} $ \ we obtain 

$$ det(\tilde{M}) = (-1)^{n+1}.(1 - \sum_{j=1}^{n+1} \ \frac{p_{j}}{r}).det(M') $$

so we conclude that  $$\sigma = \frac{r}{r - \sum_{j=1}^{n+1} \ p_{j}}.$$
 Now we have two cases :
\begin{enumerate}[i)]
\item $ r - \sum_{p_{j} < 0} \ p_{j} = d+h > \sum_{p_{j} > 0} \ p_{j} = d$. Then \ $r - \sum_{j=1}^{n+1} \ p_{j} =  (d + h) -  d = h $.

So \ $\sigma = r/h$, and the exponent of \ $\lambda$ \ in \ $P$ \ is \ $r$.
\item $ \sum_{p_{j} > 0} \ p_{j} = d + h >  r - \sum_{p_{j} < 0} \ p_{j} = d $. Then \ $r - \sum_{j=1}^{n+1} \ p_{j} = d - (d + h) = -h$.

 So \ $\sigma = - r/h$, and the exponent of \ $\lambda$ \ in \ $P$ \ is \ $-r$.
\end{enumerate}
Consider now the case of a degree \ $k$ \  monomial \ $\mu \in \C[x_{0}, \dots, x_{n}]$. Then there exists again \ $\sigma',\tau'$ \ in \ $\mathbb{Q}$ \ such that \ $(\sigma'.a + \tau'.b)[\mu] = [m_{n+2}.\mu]$ \ in \ $E_{f}$. As \ $a[\mu], (\beta_{i}+1).b[\mu], i \in [0,n]$, where \ $\beta_{i}$ \ is the degree in \ $x_{i}$ \ of \ $\mu : = x^{\beta}$, are again given from the \ $[m_{j}.\mu], j\in [1,n+2]$ \ by the same Cramer system, we conclude that \ $\sigma' = \sigma$. To conclude the proof it is enough to apply the proposition \ref{connexion}.$\hfill \blacksquare$\\

Note that in the case i) above \ $P : = P_{d+h} + c.\lambda^{r}.P_{d}$ \ annihilated \ $[\mu]$ \ in \ $E_{f}$ \ and in the case ii) we have \ $P : = P_{d+h} + c.\lambda^{-r}.P_{d}$.\\

The lemma \ref{tres facile} gives  that \ $\lambda.\nabla(P.[\mu]) = -(\sigma.a + (\tau' - k.\sigma).b).P[\mu]$ \ which makes explicit the fact that \ $\lambda.\nabla$ \ is well defined on \ $\Ab[\lambda].[\mu] \subset E_{f}$.\\

\parag{Remark} Recall that in \ $[B.13]$ \ we have built in an explicit way a differential equation in \ $s \in S$, depending in a very simple and concrete  way on \ $\lambda \in \C^{*}$ \ which is satisfied by \ $\varphi$. So it is easy to see that the knowledge of a formal asymptotic expansion when \ $s $ \ goes to \ $0$ \ in \ $S$\footnote{This is always the case when \ $S$ \ contains an open sector with edge at the origin.} \ for a given \ $\lambda_{0}$, of the type
$$ \varphi(\lambda_{0}, s) \simeq \sum_{i,j} C_{i,j}.s^{\rho_{i}}.(Log s)^{j} $$
where \ $ \ \rho_{1}, \dots,\rho_{I} \ {\rm are \ in}\  -1+ \mathbb{Q}^{*+} , \ j\in [0,n] \ {\rm are \ integers}$ \ and \ $C_{i,j}$ \ are in \ $\C[[s]]$,  determines (uniquely) via the partial differential equation above, a formal expansion of the same type for each given \ $\lambda \in \Omega$, whose coefficients  \ $C_{i,j}^{\lambda}$ \ are polynomials in \ $Log \lambda$ \ easily computable from the coefficients \ $C_{i,j}^{\lambda_{0}} : = C_{i,j}$ \ of the asymptotic expansion at the initial value  \ $\lambda_{0}$ \ of \ $\lambda$. This is described in the following lemma.\\

\begin{lemma}\label{asymptotic}
Let \ $\Omega$ \ be a simply connected domain in \ $\C^{*}$. Let \ $(\rho_{i})_{i \in I}$ \ be a finite collection of rational numbers strictly bigger than \ $-1$. Assume that the formal power serie
$$ \varphi_{\lambda} : = \sum_{k =0}^{N}\sum_{i \in I}\sum_{m \geq 0} \  c_{m}^{i,k}(\lambda).s^{m+\rho_{i}}.(Log s)^{k}\big/k! $$
where \ $c_{m}^{i,k}$ \ are holomorphic functions in \ $\Omega$, satisfies the partial differential equation
$$ \lambda.\frac{\partial }{\partial \lambda}\frac{\partial }{\partial s} \varphi_{\lambda} = \alpha.s\frac{\partial (\varphi_{\lambda})}{\partial s} + \beta.\varphi_{\lambda} $$
for each \ $\lambda \in \Omega$. Then for each \ $i,k$ \ fixed, the function \ $c_{m}^{i,k}$ \ is a polynomial in \ $Log\lambda$ \ of degree \ $\leq m$ \ for each \ $m$. Moreover the collection of  numbers \ $c_{m}^{i,k}(\lambda_{0})$ \ for a given \ $\lambda_{0} \in \Omega$ \ determines uniquely these polynomials.
\end{lemma}

\parag{Proof} The partial differential equation implies the following recursion relation for each \ $i,k,m$ :
$$ (m +\rho_{i}+ 1).\lambda.\frac{\partial c_{m+1}^{i,k}(\lambda) }{\partial \lambda} + \lambda.\frac{\partial c_{m+1}^{i,k+1}(\lambda) }{\partial \lambda} = \big(\alpha.(m + \rho_{i}) + \beta\big).c_{m}^{i,k}(\lambda) + \alpha. c_{m}^{i,k+1}(\lambda) $$
We shall make a descending induction on \ $k$. For \ $k = N$ \ the recursion relation reduces to
$$  (m +\rho_{i}+ 1).\lambda.\frac{\partial c_{m+1}^{i,N}(\lambda) }{\partial \lambda} =  \big(\alpha.(m + \rho_{i}) + \beta\big).c_{m}^{i,N}(\lambda) $$
and an easy induction on \ $m \geq 0$ \ gives our assertion for \ $k = N$.\\
Assuming the statement proved for \ $k+1$ \ a simple quadrature in \ $\lambda$ \ implies the case \ $k$.$\hfill \blacksquare$

\section{Two families of examples with \ $d = 2$ \ and \ $h = 1$.}

\subsection{The family \ $x^{2u} + y^{2v} + z^{2w} + \lambda.x^{u}.y^{v}z^{w}$.}

The condition to be in our situation is \ $u.v.w  > 0$. Then we have the relation \ $ m_{4}^{2} = \lambda^{2}.m_{1}.m_{2}.m_{3}$ \ and it shows that \ $d = 2$ \ and \ $h = 1$.\\
 Note that the only singularity of \ $f$ \ in  \ $\{ f = 0 \}$ \ is the origine.\\
To compute \ $P : = P_{3} + c.\lambda^{-2}.P_{2}$ \ which annihilates \ $[1]$ \ is not difficult. We find 
\begin{align*}
& P =  (a - (2 + \frac{u+v}{2u.v}).b)(a - (1 + \frac{u+w}{2u.w}).b))(a - (\frac{v+w}{2v.w}).b) \ + \\
& \quad  - 4\lambda^{-2}.(a - (\frac{3}{2} + \frac{u.v+v.w+w.u}{2u.v.w}).b))(a - \frac{u.v+v.w+w.u}{2u.v.w}).b)).
\end{align*}
In this case we have 
 $$ \lambda.\nabla([1]) = 2.(a - (\frac{u.v + v.w + w.u}{2u.v.w}).b)[1] = -m_{4} .$$
 Here we are in the case ii) above (so \ $\sigma = -2$).\\
Let me illustrate this family on a simple example : \ $f = x^{4} + y^{4} + z^{2} + \lambda.x^{2}.y^{2}.z $ \ corresponding to \ $u = v = 2, w = 1$. In this case we find
$$ P : = (a - \frac{5}{2}.b)\big[(a - \frac{7}{4}.b)(a - \frac{3}{4}.b) - 4.\lambda^{-2}.(a-b)\big]  \quad {\rm and} \quad \lambda.\nabla([1]) =  2(a - b)[1] .$$

\bigskip

\subsection{The family \ $ x^{2p}.z^{u} + y^{2q}.z^{v} + z^{u+v} + \lambda.x^{p}.y^{q}$.}

The condition to be in our situation is \ $p.q.(u+v) > 0 $. Note that the singularity at the origine  is not isolated in general in  these cases. We have here the equality
$$ 2.\alpha_{4} = \alpha_{1} + \alpha_{2} - \alpha_{3}$$
The relation which determines \ $P$ \ annihilating \ $[1]$ \  is given by \ $m_{4}^{2}.m_{3} = \lambda^{2}.m_{1}.m_{2}$, so \ $r = 2, d= 2, h = 1$ \ and we are in the case i).\\
The computation of \ $P$ \ gives
\begin{align*}
& P = (a - (2 +\frac{p + q}{2p.q}).b)(a - (\frac{1}{2} + \frac{p.u+q.v+2p.q}{2p.q.(u+v)}).b)(a - (\frac{p.u+q.v+2p.q}{2p.q.(u+v)}).b) \  +  \\
& \quad  -4\lambda^{2}.(a - (1 + \frac{p.u+q.v+2p.q+p.(u+v)}{2p.q.(u+v)}).b)(a - (\frac{p.u+q.v+2p.q+ q.(u+v)}{2p.q.(u+v)}).b).
\end{align*}
The computation of \ $(\sigma, \tau)$ \ such that \ $(\sigma.a + \tau.b)[1] = m_{4}$ \ is easy and it gives 
$$ \lambda.\nabla([1]) =  -  (2.a - \frac{p.u + q.v + 2p.q}{p.q.(u+v)}.b) = -m_{4} .$$
 
Again let me illustrated by an example : for \ $p = q = 2$ \ and \ $u = v = 1$ \ so for \ $f = x^{4}.z + y^{4}.z + z^{2} + \lambda.x^{2}.y^{2}$.
We find
$$ P : =  (a - \frac{5}{2}.b)(a - \frac{5}{4}.b)(a - \frac{3}{4}.b) - 4\lambda^{2}.(a - 2b)(a - b)  \quad {\rm and} \quad \lambda.\nabla([1]) = - 2(a -\frac{3}{4}.b)[1]$$

\section{Appendix}

It is interesting to remark that the proposition \ref{connexion} is a special case in a specific algebraic setting of a general result on the the filtered Gauss-Manin connexion of a holomorphic function depending holomorphically of a  parameter. This is the goal of this appendix to precise this point.\\

Let \ $M$ \ be a complex manifold, \ $D$ \ an open disc in \ $\C$ \  and let \ $f : D\times M \to \C$ \ be a holomorphic function. Denote \ $K_{/}^{p} : = Ker\big[ d_{/}f\wedge ; \Omega_{/}^{p} \to \Omega_{/}^{p+1}\big]$ \ for \ $p \ge2$ \ and \ $K_{/}^{1} : = Ker \big[ d_{/}f\wedge : \Omega_{/}^{1} \to \Omega_{/}^{2}\big]\big/\mathcal{O}.d_{/}f$ \ where \ $d_{/}f$ \ is the \ $\lambda-$relative differential of \ $f$ \ and  \ $\Omega_{/}^{p}$ \ the sheaf of \ $\lambda-$relative  holomorphic \ $p-$forms (compare with [B.08]).\\

Denote  by \ $(K_{/}^{\bullet},d_{/})$ \ the topological  restriction  of the \ $\lambda-$relative de Rham  complex (defined above) for the map \ $(\lambda,x) \mapsto (\lambda, f(\lambda,x))$, to the analytic subset 
$$Z : = \{ d_{/}f = 0\} $$  
 and let \ $\mathcal{H}^{p}$ \ the \ $p-$th cohomology sheaf of this complex. Recall that these cohomogy sheaves have a natural structure of left \ $\Ab[\lambda]-$modules with the action of \ $a$ \  given by the multiplication by \ $f$ \   and with  the action of \ $b$ \  defined by \ $d_{/}f\wedge d_{/}^{-1}$.

\begin{prop}\label{general}
There exists a natural graded map \ $\nabla^{\bullet} : \mathcal{H}^{\bullet} \to \mathcal{H}^{\bullet}$ \ with the following properties :
\begin{enumerate}
\item For \ $\omega = d_{/}\xi \in K_{/}^{p+1}\cap Ker\ d_{/}$ \ we have \ $\nabla([\omega]) = [d_{/}f\wedge\frac{\partial \xi}{\partial \lambda} -\frac{\partial f}{\partial \lambda}.\omega] .$
\item The map \ $b^{-1}.\nabla$ \ well defined on \ $ b.\tilde{\mathcal{H}}^{\bullet}$ \ where \ $\tilde{\mathcal{H}}^{\bullet}\ : = \mathcal{H}^{\bullet}  \big/(b-torsion)$, with value in \ $\tilde{\mathcal{H}}^{\bullet}$, commutes with \ $a$ \ and \ $b$ \ and is a \ $\lambda-$connection.
\item If \ $(\gamma_{s,\lambda})_{(s,\lambda) \in S\times \Omega}$ \ is a horizontal family of  compact\ $p-$cycles in the fibers of \ $(f,\lambda)$ \ over an open set in \ $D\times M \setminus C(f,\lambda)$ \ where \ $C(f,\lambda)$ \ is the set of critical values of the map \ $(\id, f ) : D\times M \to D \times \C $, we have for any \ $\omega \in  K_{/}^{p+1}\cap Ker\ d_{/}$ \ the equality
$$  \frac{\partial }{\partial s} \frac{\partial }{\partial \lambda}\big[ \int_{\gamma_{s,\lambda}} \ \frac{\omega}{d_{/}f } \big] =  \int_{\gamma_{s,\lambda}}  \ \frac{ \nabla(\omega)}{d_{/}f }.$$
\end{enumerate}
\end{prop}

\parag{Proof} First remark that if \ $\omega = d_{/}\xi \in K^{p+1}$ \  with \ $d_{/}f\wedge \xi = 0$, we have \ $d_{/}(\frac{\partial f}{\partial \lambda})\wedge \xi + d_{/}f\wedge\frac{\partial \xi}{\partial \lambda} = 0$ \ so
 \ $ \nabla(d_{/}\xi) = - d_{/}(\frac{\partial f}{\partial \lambda}.\xi) $ \ is in \ $d_{/}K^{p}$.\\
 Now for  \ $\omega = d_{/}\xi \in K^{p+1}$ \  we have \ $d_{/}f\wedge \nabla(d_{/}\xi) = 0$ \ and 
  $$d_{/}\big(\nabla(d_{/}\xi)\big) =  -d_{/}f\wedge d_{/}\big(\frac{\partial \xi}{\partial \lambda} \big) - d_{/}\big( \frac{\partial f }{\partial \lambda}\big)\wedge d_{/}\xi $$
  and we obtain that
  $$d_{/}\big(\nabla(d_{/}\xi)\big) =  - \frac{\partial }{\partial \lambda}\big(d_{/}f \wedge \omega\big) = 0 $$
  using the fact that \ $ \frac{\partial }{\partial \lambda}\big(d_{/}f \wedge \omega\big) \equiv 0$.\\
  The proof of the other statements are analogous to the corresponding ones in proposition \ref{connexion}. $\hfill \blacksquare$
  
  \parag{Remarks} \begin{enumerate}
  \item As above we have a more precise formulation for the properties in assertion 2. of the proposition above with the  following relations in \ $\mathcal{H}^{\bullet}$
  \begin{align*}
  &  \nabla(\varphi.\omega) = \frac{\partial \varphi}{\partial \lambda}.b.\omega + \varphi.\nabla(\omega) \quad {\rm for} \quad \varphi \in \mathcal{O}_{\lambda} \quad {\rm and} \quad \omega \in \mathcal{H}^{\bullet}\\
  & b.\nabla = \nabla.b \quad{\rm and} \quad  \nabla.a = (a - b).\nabla.
  \end{align*}
  \item The generalization of this proposition to several holomorphic parameters is immediate.$\hfill \square$
  \end{enumerate}

  \parag{Bibliography}\begin{itemize}
  \item{[B.08]}  Barlet,D.  \textit{Sur certaines singularit\'es d'hypersurfaces II},  \\ Journal of Algebraic 
  Geometry 17 (2008), p.199-254.
  \item{[B.09]} Barlet,D. {\it P\'eriodes \'evanescentes et (a,b)-modules monog\`enes}, Bollettino U.M.I. (9) II (2009) p.651-697.
  \item{ [B.13]} Barlet, D. {\it Algebraic differential equations associated to some \\ polynomials}  arXiv:1305.6778 ( math.arxiv 2013).
  \end{itemize}

  \end{document}